\documentclass[10pt,reqno]{amsart}

\usepackage{amsthm}

\usepackage{eucal} 

\usepackage{enumerate}
\usepackage[longnamesfirst,numbers]{natbib}
\usepackage[bookmarks,colorlinks=true,citecolor=blue]{hyperref}

\theoremstyle{plain}
\newtheorem{theorem}{Theorem}[section]
\newtheorem{lemma}[theorem]{Lemma}

\newtheorem{proposition}[theorem]{Proposition}

\theoremstyle{definition}
\newtheorem{definition}[theorem]{Definition}

\theoremstyle{remark}

\numberwithin{equation}{section}

\title{Probability Measures and Effective Randomness}

\author{Jan Reimann}

\address{Department of Mathematics  \\
  University of California, Berkeley \\
}

\email{reimann@math.berkeley.edu}

\author{Theodore A. Slaman} 

\address{Department of Mathematics  \\
  University of California, Berkeley \\
}

\email{slaman@math.berkeley.edu}

\date{}

\newcommand{\Nat}{\ensuremath{\mathbb{N}}}

\newcommand{\Rat}{\ensuremath{\mathbb{Q}}}

\newcommand{\Cant}{\ensuremath{2^{\omega}}}

\newcommand{\Str}[1][<\omega]{\ensuremath{2^{#1}}}

\newcommand{\Sle}{\ensuremath{\subset}}
\newcommand{\Sleq}{\ensuremath{\subseteq}}

\newcommand{\Sgeq}{\ensuremath{\supseteq}}

\newcommand{\Cl}[1]{\ensuremath{#1}}

\newcommand{\Cyl}[1]{\ensuremath{N_{#1}}}
\newcommand{\ACyl}[1]{\ensuremath{N(#1)}}

\newcommand{\Fam}[1]{\ensuremath{\mathcal{#1}}}

\newcommand{\Rest}[1]{\ensuremath{\, \lceil #1}}
\newcommand{\Op}[1]{\ensuremath{\operatorname{#1}}}

\newcommand{\Conc}{\ensuremath{\mbox{}^\frown}}
\newcommand{\Estr}{\ensuremath{\epsilon}}
\newcommand{\Tup}[1]{\ensuremath{\langle #1 \rangle}}

\newcommand{\join}[1][\mbox{}]{\ensuremath{\oplus_{#1}}}

\newcommand{\Pmeas}{\ensuremath{\CMcal{P}}}
\newcommand{\Leb}{\ensuremath{\mathcal{L}}}
\newcommand{\ZFC}{\ensuremath{\mathsf{ZFC}}}

\DeclareMathOperator{\NCR}{NCR}
\DeclareMathOperator{\T}{T}
\DeclareMathOperator{\TT}{tt}

\begin{document}

\maketitle

\begin{abstract}	
  We study the question, ``For which reals $x$ does there exist a
  measure $\mu$ such that $x$ is random relative to $\mu$?''  We show
  that for every nonrecursive  $x$, there is a measure which makes
  $x$ random without concentrating on $x$.  We give several conditions
  on $x$ equivalent to there being continuous measure which makes $x$
  random.  We show that for all but countably many reals $x$ these
  conditions apply, so there is a continuous measure which makes $x$
  random.  There is a meta-mathematical aspect of this investigation.
  As one requires higher arithmetic levels in the degree of
  randomness, one must make use of more iterates of the power set of
  the continuum to show that for all but countably many $x$'s there is
  a continuous $\mu$ which makes $x$ random to that degree.
\end{abstract}

%
%
\section{Introduction} \label{sec:intro}
Most studies on algorithmic randomness focus on reals random with
respect to the uniform distribution, i.e.\ the $(1/2,1/2)$-Bernoulli
measure, which is measure theoretically isomorphic to Lebesgue measure
on the unit interval. The theory of uniform randomness, with all its
ramifications (e.g.\ computable or Schnorr randomness) has been well
studied over the past decades and has led to an impressive theory.

Recently, a lot of attention focused on the interaction of algorithmic
randomness with recursion theory: What are the computational
properties of random reals? In other words, which computational
properties hold effectively for almost every real? This has led to a
number of interesting results, many of which will be covered in a
forthcoming book by \citet{downey-hirschfeldt:ip}.

While the understanding of ``holds effectively'' varied in these
results (depending on the underlying notion of randomness, such as
computable, Schnorr, or weak randomness, or various arithmetic levels
of Martin-L\"of randomness, to name only a few), the meaning of ``for
almost every'' was usually understood with respect to Lebesgue
measure.  One reason for this can surely be seen in the fundamental
relation between uniform Martin-L\"of tests and descriptive complexity
in terms of (prefix-free) Kolmogorov complexity: A real is not covered
by any Martin-L\"of test (with respect to the uniform distribution) if
and only if all of its initial segments are incompressible (up to a
constant additive factor).

However, one may ask what happens if one changes the underlying
measure. This question is virtually as old as the theory of
randomness. \citet{martinloef:1966} defined randomness not only for
Lebesgue measure but also for arbitrary Bernoulli
distributions. Levin's contributions in the 1970's
\citep{zvonkin-levin:1970, levin:1973,levin:1974, levin:1976} extended
this to arbitrary probability measures. He obtained a
number of remarkable results and principles such as the existence of
uniform tests, conservation of randomness, and the existence of
neutral measures. 

In this paper we will survey a recent line of research by the authors which dealt with the question \emph{for which reals $x$ does there exist a probability measure which makes $x$ random without concentrating on $x$}. We consider two kinds measures -- arbitrary probability measures, which may have atoms (reals other than $x$ on which the measure concentrates), and \emph{continuous measures}, i.e. non-atomic measures. The investigations exhibit an interesting, and quite unexpected, connection between the randomness properties of a real and its logical complexity, in the sense of recursion or set theoretic hierarchies. In the following we will try to describe this connection in some detail. We will sketch proofs to provide some intuition, but for a full account we have to refer the reader to the forthcoming research papers \citep{reimann-slaman:ip2,reimann-slaman:tbs}.

%
%
\section{Measures and Randomness} \label{sec-meas-cant}

In this section we introduce the basic notions of measure on the
Cantor space $\Cant$ and define randomness for arbitrary probability
measures.

The \emph{Cantor space} $\Cant$ is the set of all infinite binary
sequences, also called \emph{reals}. The topology generated by the
\emph{cylinder sets}
\[
N_\sigma = \{ x : \: x\Rest{n} = \sigma\},
\]
where $\sigma$ is a finite binary sequence, turns $\Cant$ into a
compact Polish space. We will occasionally use the notation
$N(\sigma)$ in place of $\Cyl{\sigma}$ to avoid multiple
subscripts. $\Str$ denotes the set of all finite binary sequences. If
$\sigma, \tau \in \Str$, we use $\Sleq$ to denote the usual prefix
partial ordering. This extends in a natural way to $\Str \cup
\Cant$. Thus, $x \in \Cyl{\sigma}$ if and only if $\sigma \Sle
x$. Finally, given $U \subseteq \Str$, we write $\Cyl{U}$ to denote
the open set induced by $U$, i.e. $\Cyl{U} = \bigcup_{\sigma \in U}
\Cyl{\sigma}$.

%
%
\subsection{Probability measures}

A \emph{probability measure} on $\Cant$ is a \emph{countably additive,
  monotone function} $\mu: \Fam{F} \to [0,1]$, where $\Fam{F}
\subseteq \mathcal{P}(\Cant)$ is $\sigma$-algebra and $\mu(\Cant) =
1$. $\mu$ is called a \emph{Borel probability measure} if $\Fam{F}$ is
the Borel $\sigma$-algebra of $\Cant$. It is a basic result of measure
theory that a probability measure is uniquely determined by the values
it takes on an algebra $\Fam{A} \subseteq \Fam{F}$ that generates
$\Fam{F}$. It is not hard to see that the Borel sets are generated by
the algebra of \emph{clopen sets}, i.e.\ finite unions of basic open
cylinders. Normalized, monotone, countably additive set functions on
the algebra of clopen sets are induced by any function $\rho: \Str \to
[0,1]$ satisfying
\begin{equation} \label{equ-probability measure}
  \rho(\Estr) = 1 \qquad \text{and} \qquad \rho(\sigma) = \rho(\sigma \Conc 0) + \rho(\sigma \Conc 1)
\end{equation}
for all finite sequences $\sigma$. Then $\mu(\Cyl{\sigma}) =
\rho(\sigma)$ yields an monotone, additive function on the clopen
sets, which in turn uniquely extends to a Borel probability measure on
$\Cant$. In the following, we will deal exclusively with Borel
probability measures, and hence we will identify such measures with
the underlying function on cylinders satisfying \eqref{equ-probability
  measure}, and write, in slight abuse of notation, $\mu(\sigma)$
instead of $\mu(\Cyl{\sigma})$. Besides, we will mostly speak of
\emph{measures}, understanding Borel probability measures.

The \emph{Lebesgue measure} $\Leb$ on $\Cant$ is obtained by
distributing a unit mass uniformly along the paths of $\Cant$, i.e.\
by setting $\Leb(\sigma) = 2^{-|\sigma|}$. A \emph{Dirac measure}, on
the other hand, is defined by putting a unit mass on a single real,
i.e.  for $x \in \Cant$, let
\[
\delta_x(\sigma) = \begin{cases}
  1 & \text{if } \sigma \subset x, \\
  0 & \text{otherwise.}
\end{cases}
\]
If, for a measure $\mu$ and $x \in \Cant$, $\mu(\{x\}) > 0$, then $x$
is called an \emph{atom} of $\mu$. Obviously, $x$ is an atom of
$\delta_x$. A measure that does not have any atoms is called
\emph{continuous}.

%
%
\section{Martin-L\"of Randomness} \label{sec-randomness}

It was Martin-L\"of's fundamental idea to define randomness by
choosing a \emph{countable family} of nullsets. For any non-trivial
measure, the complement of the union of these sets will have positive
measure, and any point in this set will be considered
\emph{random}. There are of course many possible ways to pick a
countable family of nullsets. In this regard, it is very benefiting to
use the framework of recursion theory and effective descriptive set
theory.

%
%
\subsection{Nullsets}

Before we go on to define Martin-L\"of randomness formally, we note
that every nullset is contained in a $G_\delta$-nullset.

\begin{proposition}
  Suppose $\mu$ is a measure. Then a set $A \subseteq \Cant$ is
  $\mu$-null if and only if there exists a set $U \subseteq \Nat
  \times \Str$ such that for all $n$,
  \begin{equation} \label{equ-nullset}
    A \subseteq \ACyl{U_n} \quad \text{and} \quad 
    \sum_{\sigma \in U_n} \mu(\Cyl{\sigma}) \leq 2^{-n},
  \end{equation} 
  where $U_n = \{ \sigma : \: (n, \sigma) \in U \}$.
\end{proposition}

Of course, the $G_\delta$-cover of $A$ is given by $\bigcap_n U_n$. 

%
%
\subsection{Martin-L\"of tests and randomness}

Essentially, a Martin-L\"of test is an effectively presented
$G_\delta$ nullset (relative to some parameter $z$).

\begin{definition} \label{def-test} Suppose $z \in \Cant$ is a real. A
  \emph{test relative to $z$}, or simply a \emph{$z$-test}, is a set
  $W \subseteq \Nat \times \Str$ which is recursively enumerable in
  $z$. Given a natural number $n \geq 1$, an \emph{$n$-test} is a test
  which r.e.\ in $\emptyset^{(n-1)}$, the $(n-1)$st Turing jump of the
  empty set. A real $x$ \emph{passes} a test $W$ if $x \not\in
  \bigcap_n \ACyl{W_n}$.
\end{definition}

Passing a test $W$ means not being contained in the $G_\delta$ set
given by $W$. The condition `\emph{r.e.\ in $z$}' implies that the
open sets given by the sets $W_n$ form a uniform sequence of
$\Sigma^0_1(z)$ sets, and the set $\bigcap_n \ACyl{W_n}$ is a
$\Pi^0_2(z)$ subset of $\Cant$.

To test for randomness with respect to a measure, we have to ensure
two things: First that a test $W$ actually describes a
nullset. Second, that the information present in a measure is
available to the test. Te first criterion we call \emph{correctness}.

\begin{definition} 
  Suppose $\mu$ is a measure on $\Cant$. A test $W$ is \emph{correct
    for $\mu$} if  
  \begin{equation} \label{equ-correct-test}
    \sum_{\sigma \in W_n} \mu(\Cyl{\sigma}) \leq 2^{-n}.
  \end{equation}
\end{definition}

To incorporate measures into an effective test for randomness we have
to represent it in a form that makes it accessible for recursion
theoretic methods. Essentially, this means to code a measure via an
infinite binary sequence or a function $f:\Nat \to
\Nat$. Unfortunately, there are many possible such
representations. Hence, strictly speaking, we will deal with
\emph{randomness with respect to a representation} of a measure, not
the measure itself. However, we will see that for one of our main
topics, randomness for continuous measures, representational issues
can be resolved quite elegantly.

The most straightforward representation of a measure is the following.  

\begin{definition}\label{def-rational-repr}
  Given a measure $\mu$, define its \emph{rational representation}
  $r_\rho$ by letting, for all $\sigma \in \Str$, $q_1, q_2 \in \Rat$,
  \begin{equation}
    \Tup{\sigma, q_1, q_2} \in r_\rho \;
    \Leftrightarrow \; q_1 < \rho(\sigma) < q_2. 
  \end{equation}
\end{definition}

The rational representation does not reflect the topological
properties of the space of probability measures on $\Cant$. The space of probability measures $\Pmeas$ on $\Cant$ is a compact
polish space (see \citet{parthasarathy:1967}). The topology is the \emph{weak topology}, which can be
metrized by the \emph{Prokhorov metric}, for instance. There is an
\emph{effective dense subset}, given as follows: Let $Q$ be the set of
all reals of the form $\sigma\Conc 0^\omega$. Given $\bar{q} = (q_1,
\dots, q_n) \in Q^{< \omega}$ and non-negative rational numbers
$\alpha_1, \dots, \alpha_n$ such that $\sum \alpha_i = 1$, let
\[
\delta_{\bar{q}} = \sum_{k=1}^n \alpha_k \delta_{q_k},
\]
where $\delta_x$ denotes the \emph{Dirac point measure} for $x$. Then
the set of measures of the form $\delta_{\bar{q}}$ is dense in
$\Pmeas$.

The recursive dense subset $\{\delta_{\bar{q}}\}$ and the
effectiveness of the metric $d$ between measures of the form
$\delta_{\bar{q}}$ suggests that the representation reflects the
topology effectively, i.e.\ the set of representations should be
$\Pi^0_1$. However, this is not true for the set of rational
representations of probability measures. Instead, we have to resort to
other representations in metric spaces, such as Cauchy
sequences. Using the framework of \emph{effective descriptive set
  theory}, as for example presented in \citet{moschovakis:1980}, one
can obtain the following.

\begin{theorem} \label{thm-pi01-repres-meas}
  There is a recursive surjection
  \[
  \pi: \: \Cant \to \Pmeas
  \]
  and a $\Pi^0_1$ subset $P$ of $\Cant$ such that $\pi\Rest{P}$ is
  one-one and $\pi(P) = \Pmeas$.
\end{theorem}

In the following sections, we will always assume that a measure is
either represented by its rational representation or via the the set
$P \subseteq \Cant$ of the previous theorem. The definition of
randomness, however, works for any representation.

\begin{definition}\label{def-test-relative-to-measure}
  Suppose $\mu$ is a probability measure on $\Cant$, $\rho_\mu \in \Cant$ is a
  representation of $\mu$, and $z \in \Cant$ is a real.  A real is
  \emph{Martin-L\"of $n$-random for $\mu$ relative to $\rho_\mu$ and
    $z$}, or simply \emph{$(n,z)$-random for $\rho_\mu$} if it passes
  all $(\rho_\mu \join z)^{(n-1)}$-tests which are correct for $\mu$.
\end{definition}

If the representation is clear from the context, we speak of
$(n,z)$-random\-ness for $\mu$. If $\mu$ is Lebesgue measure $\Leb$, we
drop reference to the measure and simply say ``$x$ is
$(n,z)$-random''. We also drop the index $1$ in case of
$(1,z)$-randomness and simply speak of ``randomness relative to $z$''
or $z$-randomness.

\medskip
Since there are only countably many Martin-L\"of $n$-tests, it follows
from countable additivity  that the set of Martin-L\"of $n$-random reals for $\mu$
has $\mu$-measure $1$. Hence there always exist $(n,z)$-random reals for any measure
$\mu$.

%
%
\subsection{Image measures and conservation of
  randomness} \label{subsec-trans-meas}

One can obtain new measures from given measures by transforming them
with respect to a sufficiently regular function. Let $f: \Cant \to
\Cant$ be a Borel (measurable) function, i.e.\  for every Borel set
$A$, $f^{-1}(A)$ is Borel, too. If $\mu$ is a measure on $\Cant$ and
$f$ is Borel, then the \emph{image measure} $\mu_f$ is defined by
\[
\mu_f(A) = \mu(f^{-1}(A)).
\] 

It can be shown that every probability measure can be obtained from
Lebesgue measure $\Leb$ by means of a measurable transformation.

\begin{theorem}[folklore, see e.g.\
  \citet{billingsley:1995}] \label{thm-transf-meas}
  If $\mu$ is a Borel probability measure on $\Cant$, then there
  exists a measurable $f: \Cant \to \Cant$ such that $\mu =
  \Leb_{f}$.
\end{theorem}

If the transformation of $\Leb$ is effective, then $f$ maps an
$\Leb$-random real to a $\Leb_{f}$-random real. This principle is
called \emph{conservation of randomness}, first introduced by
Levin. We can use it to construct measures for which a given real is
random, as we will see in the next sections.

%
%
\section{Randomness of Non-Recursive Reals}

If $x$ is an atom of some probability measure $\mu$, it is trivially
$\mu$-random. Interestingly, the recursive reals are exactly those for
which this is the only way to become random.

\begin{theorem}[\citet{reimann-slaman:tbs}] \label{thm-nonrec-rand}
  For any real $x$, the following are equivalent.
  \begin{enumerate}[(i)]
  \item There exists a (representation of a) probability measure $\mu$
    such that $\mu(\{x\}) = 0$ and $x$ is $\mu$-random.
  \item  $x$ is not recursive.
  \end{enumerate}
\end{theorem}

\begin{proof}[Proof sketch]
  If $x$ is recursive and $\mu$ is a measure with $\mu(\{x\}) = 0$, then
  we can obviously construct a $\mu$-test that covers $x$, by
  computing (recursively in $\mu$) the measure of initial segments of
  $x$, which tends to $0$.

  Now assume $x$ is not recursive. 
  A fundamental result by \citet{kucera:1985} ensures that every
  Turing degree above $\emptyset'$ contains a $\Leb$-random real. This
  result relativizes. Hence one can combine it with the
  \emph{Posner-Robinson Theorem} \citep{posner-robinson:1981}, which
  says that for every non-recursive real $x$ there exists a $z$ such
  that $x \oplus z =_{\T} z'$. This way we  obtain a real $R$ which is
  \begin{enumerate}[(1)]
  \item $\Leb$-random relative to some $z \in \Cant$, and 
  \item $\T(z)$-equivalent to $x$. 
  \end{enumerate}
  There are Turing functionals $\Phi$ and $\Psi$ recursive in $z$ such
  that
  \begin{displaymath}
    \Phi(R) = x \quad \text{ and } \quad \Psi(x) = R.	
  \end{displaymath}	
  We can use the functionals to define a class of measures that are
  possible candidates to
  render $x$ random. Given $\sigma \in \Cant$, define the set
  $\Op{Pre}(\sigma)$ to be
  the set of minimal elements of
  \begin{equation*}
    \{ \tau \in \Str: \: \Phi(\tau) \Sgeq \sigma \: \text{ and } \:
    \Psi(\sigma) \Sleq \tau \}.
  \end{equation*}
  We define a set of measures $M$ by requiring that $\mu \in M$ if and
  only if
  \begin{equation}\label{entrop:equ_meas-cond}
    \forall\sigma[\Leb(\Op{Pre}(\sigma)) \leq \mu(\sigma) \leq
    \Leb(\Psi(\sigma))]. 
  \end{equation}
  The first inequality ensures that $\mu$ \emph{dominates} an image
  measure induced
  by $\Phi$. This will ensure that any Martin-L\"of random real is
  mapped by $\Phi$ to a $\mu$-random real. The second inequality
  guarantees that $\mu$ is non-atomic on the domain of $\Psi$.

  One can show the topological representations of the measures in $M$
  (Theorem \ref{thm-pi01-repres-meas}) form a non-empty $\Pi^0_1$
  class $\Cl{M}$ in $\Cant$
  relative to $z$. 

  In order to apply conservation of randomness, we have to know that
  one of the measures in $M$, when given as an additional information
  to a $\Leb$-test, will not destroy the randomness of $R$. This is
  ensured by the following basis result for $\Pi^0_1$ sets regarding
  relative randomness (essentially a consequence of
  \emph{compactness}).
\end{proof}	

\begin{theorem}[\citet{reimann-slaman:tbs},
  \citet{downey-hirschfeldt-miller-nies:2005}] \label{thm-basis-rel-rand} 
  Let $S$ be $\Pi^0_1(z)$. If $R$ is $\Leb$-random relative to
  $z$, then there exists $y \in S$ such that $R$ is $\Leb$-random
  relative to $y \oplus z$.
\end{theorem}

%
%
\section{Randomness for continuous measures}

A natural question arising in the context Theorem
\ref{thm-nonrec-rand} is whether the measure making a real random can
be ensured to have certain regularity properties; in particular, can
it be chosen \emph{continuous}? 

\citet{reimann-slaman:tbs} gave an explicit construction of a
non-recursive real not random with respect to any continuous
measure. Call such reals \emph{$1$-ncr}. In general, let $\NCR_n$ be
the set of reals which are not $n$-random with respect to any
continuous measure.

\citet{kjoshanssen-montalban:2005} observed that any member of a
countable $\Pi^0_1$ class is an element of $\NCR_1$. 

\begin{proposition}
  If $A \subseteq \Cant$ is $\Pi^0_1$ and countable, then no member of
  $A$ can be in $\NCR_1$.
\end{proposition}

\begin{proof}[Proof idea]
  If $\mu$ is a continuous measure, then obviously $\mu(A) = 0$. One
  can use a recursive tree $T$ such that $[T] = A$ to obtain a
  $\mu$-test for $A$. 
\end{proof}

It follows from results of \citet{cenzer-etal:1986} that members of
$\NCR_1$ can be found throughout the hyperarithmetical hierarchy of
$\Delta^1_1$, whereas \citet{kreisel:1959} had shown earlier that each
member of a countable $\Pi^0_1$ class is in fact hyperarithmetical.

Quite surprisingly, $\Delta^1_1$ turned out to be the precise upper
bound for $\NCR_1$. An analysis of the proof of Theorem
\ref{thm-nonrec-rand} shows that if $x$ is \emph{truth-table}
equivalent to a $\Leb$-random real, then the ``pull-back'' procedure
used to devise a measure for $x$ yields a continuous measure. More
generally, we have the following.

\begin{theorem}[\citet{reimann-slaman:tbs}] \label{thm-char-cont-rand}
  Let $x$ be a real. For any $z \in \Cant$ and any $n \geq 1$, the
  following are equivalent.
  \begin{enumerate}[(i)]
  \item $x$ is $(n,z)$-random for a continuous measure $\mu$ recursive
    in $z$.
  \item $x$ is $(n,z)$-random for a continuous dyadic measure $\nu$
    recursive in $z$.
  \item There exists a functional $\Phi$ recursive in $z$ which is an
    order-preserving homeomorphism of $\Cant$ such that $\Phi(x)$ is
    $(n,z)$-random.
  \item $x$ is truth-table equivalent relative to $z$ to a
    $(n,z)$-random real.
  \end{enumerate}
\end{theorem}

Here \emph{dyadic} measure means that the values of $\mu$ on the open
cylinders are of the form $\mu(\sigma) = m/2^n$ with $m,n \in
\Nat$. The theorem can be seen as an effective version of the
\emph{classical isomorphism theorem} for continuous probability
measures (see for instance \citet{kechris:1995}).\footnote{The theorem
  suggests that for continuous randomness representational issues do
  not really arise, since there is always a measure with a
  computationally minimal representation.} 

\citet{woodin:sub}, using a variation on Prikry forcing, was able to prove that if $x
\in \Cant$ is not hyperarithmetic, then there is a $z \in \Cant$ such
that $x\oplus z \equiv_{\TT(z)} z'$, i.e.\ outside $\Delta^1_1$ the
Posner-Robinson theorem holds with truth-table equivalence. Hence we
can infer the following result.

\begin{theorem}[\citet{reimann-slaman:tbs}] \label{thm-ncr1-delta11}
  If a real $x$ is not $\Delta^1_1$, then there exists a continuous
  measure $\mu$ such that $x$ is $\mu$-random.
\end{theorem}

It is on the other hand an open problem whether every real in $\NCR_1$
is a member of a countable $\Pi^0_1$ class.

\medskip
One may ask how the complexity and size of $\NCR_n$ grows with $n$. It
turned out all levels of $\NCR$ are countable.

\begin{theorem}[\citet{reimann-slaman:ip2}]
  For all $n$, $\NCR_n$ is countable.
\end{theorem}

\begin{proof}[Proof idea]
  The first step is to use \emph{Borel determinacy} to show that the
  complement of $\NCR_n$ contains an upper Turing cone. This follows
  from the fact that the complement of $\NCR_n$ contains a Turing
  invariant and cofinal (in the Turing degrees) Borel set, which can
  be seen as follows. 

  If for two reals $x,y$, $x \equiv_{\T(z)} y$, then $x
  \equiv_{\TT(z')} y$. Suppose $x \equiv_{\T(z)} R$ where $R$ is
  $(n+1)$-random relative to $z$. Then, since $R$ is  $n$-random
  relative to $z'$, it follows from Theorem \ref{thm-char-cont-rand}
  that $x$ is random with respect to some continuous measure.

  So if we let $B \subseteq \Cant$ be the set
  \[
  \{ x \in \Cant: \: \exists z \, \exists R \, (x \equiv_{\T} z \oplus
  R \;\; \& \;\; R \text{ is $(n+1)$-random relative to $z$})\},
  \]
  $B$ is a Turing invariant Borel set cofinal in the Turing degrees.
  It follows from Borel Determinacy \citep{martin:1968} that $B$
  contains an upper cone in the Turing degrees.

  The next step is to show that the elements of $\NCR_n$ show up at a
  \emph{countable level of the constructible universe} $L$. It holds
  that $\NCR_n \subseteq L_{\beta_n}$, where $\beta_n$ is the least
  ordinal such that 
  \[
  L_{\beta_n} \models \ZFC^{-}_n,
  \]
  where $\ZFC^{-}_n$ is $\ZFC$ with the power set axiom replaced by
  the existence of $n$ iterates of the power set of $\omega$. Note
  that $L_{\beta_n}$ is the level of constructibility capturing
  Martin's construction of a winning strategy in a
  $\Sigma^0_n$-game.

  Given $x \not\in L_{\beta_n}$, construct a set $G$ such that
  $L_{\beta_n}[G]$ is a model of $\mathsf{ZFC}^{-}_n$, and for all $y
  \in L_{\beta_n}[G] \cap \Cant$, $y \leq_{\T} x \oplus G$. $G$ is
  constructed by \emph{Kumabe-Slaman forcing} (see
  \citep{shore-slaman:1999}). This notion of forcing provides a method
  to extend the Posner-Robsinson Theorem to higher levels of the jump
  and beyond.  The existence of $G$ allows to conclude: If $x$ is not
  in $L_{\beta_n}$, it will belong to every cone with base in
  $L_{\beta_n}[G]$. In particular, it will belong to the cone given by
  the Borel Turing determinacy argument (relativized to $G$, here one
  has to use absoluteness), i.e.\ the cone avoiding $\NCR_n$.  Hence
  $x$ is random relative to $G$ for some continuous $\mu$, an thus in
  particular $\mu$-random.
\end{proof}

The proof of the countability of $\NCR_n$ makes essential use of Borel
determinacy.  It is known from a result by \citet{friedman:1970} that
the use of $\omega_1$-many iterates of the power set of $\omega$ are
necessary to prove Borel determinacy. In the simplest case, Friedman
showed that $\mathsf{ZFC}^{-}$ does not prove the statement ``All
${\Sigma^0_5}$-games on countable trees are determined.''  The proof
works by showing that there is a model of $\mathsf{ZFC}^{-}$ for which
${\Sigma^0_5}$-determinacy does not hold. This model is just
$L_{\beta_0}$.  The analysis extends to higher levels of the Borel
hierarchy, applying to more and more iterates of the power set.

The question suggests itself whether the proof of the countability of
$\NCR_n$ requires a similar set theoretic complexity.

\begin{theorem}[\citet{reimann-slaman:ip2}]
  For every $k$, the statement
  \begin{center}
    For every $n$, $\NCR_n$ is countable.
  \end{center}
  cannot be proven in $\mathsf{ZFC}^{-}_k$. 
\end{theorem}

\begin{proof}[Proof sketch] We show that for every fixed $k$, some
  $\NCR_n$ is cofinal in the Turing degrees of $L_{\beta_k}$. In fact,
  Jensen's \emph{master codes} \citep{jensen:1972} for $L$, the
  universe of constructible sets, are the cofinal set.

  $L$ is generated by transfinite recursion in which the recursion
  steps are closing under first order definability and forming unions.
  The master codes represent the initial segments of $L$ and are
  generated by iterating the Turing jump and taking $L$-least
  representations of direct limits.  In short, a master code is either
  definable relative to an earlier master code or is the code for the
  well-founded limit of structures each of which is coded by an earlier
  master code.  The number of iterates of the power set present in the
  initial segment of $L$ which is being coded is linked to the
  complexity of describing the direct limit used to form its
  master code.  For $\alpha$ less than $\beta_k$, there is a fixed
  bound on this complexity.  We let $M_\alpha$ denote the master code
  for $L_\beta$.

  Neither of these cases is consistent with randomness, as indicated
  by the following lemmas.

  \begin{lemma}\label{jump}
    Suppose that $n\geq 2$, $y\in\Cant$, and $x$ is $n$-random
    for $\mu$.  If $i < n$, $y$ is recursive in $(x \oplus
    \mu)$ and recursive in $\mu^{(i)}$, then $y$ is recursive in
    $\mu$.
  \end{lemma}

  \begin{lemma}\label{limit}
    Suppose that $x$ is $(n+5)$-random for $\mu$, $\prec$ is a linear ordering that is
    $\Delta^0_{n+1}$ relative to $\mu$, and $I$ is the largest initial
    segment of $\prec$ which is well-founded.  If $i<n$ and $I$ is
    $\Sigma^0_i$ in $(x\oplus\mu)$, then $I$ is recursive in $\mu$.
  \end{lemma}

  Suppose $\alpha$ less than $\beta_k$ and a continuous measure $\mu$
  are given so that $M_\alpha$ is random relative to $\mu$.
  Heuristically, we argue as follows.  We proceed by induction on
  $\beta\leq\alpha$ to prove that $M_\beta$ is recursive in $\mu$.  If
  $\beta$ is a successor, then $M_\beta$ is arithmetic in some earlier
  master code, with a uniform upper bound on the complexity of the
  definition depending on $k$.  Then, $M_\beta$ is uniformly arithmetic
  in $\mu$ and recursive in $M_\alpha$, Lemma~\ref{jump} applies.
  Otherwise, $M_\beta$ is the well-founded direct limit of structures
  recursive in $\mu$ and recursive in in $M_\alpha$, so
  Lemma~\ref{limit} applies.  In either case, $M_\beta$ is recursive
  in $\mu$.  By induction, $M_\alpha$ is itself recursive in $\mu$ and
  not $\mu$-random, a contradiction.
\end{proof}

\bibliographystyle{abbrvnat}

\end{document}